\documentstyle[12pt]{article}
\newcommand{\be}{\begin{equation}}
\newcommand{\ee}{\end{equation}}
\newcommand{\bea}{\begin{eqnarray}}
\newcommand{\eea}{\end{eqnarray}}
\newcommand{\beann}{\begin{eqnarray*}}
\newcommand{\eeann}{\end{eqnarray*}}

\newcommand{\norm}[1]{\mbox{$\|#1\|$}}

\begin{document}

\title{Rings With Topologies Induced by Spaces of Functions}
\author{R{\u{a}}zvan Gelca}
\maketitle

\hspace {15 mm}{\bf Abstract:} By considering topologies on Noetherian 
rings that carry the properties of those induced by spaces of functions,
we prove that if an ideal is closed then every prime ideal associated
to it is closed (thus answering a question raised in [5]). The converse
is also true if we assume that a topological version of the Nullstellensatz
holds, and we prove such a result for the ring of polynomials in two variables 
endowed with the topology induced by the Hardy space.
 The topological completion of the ring is a module, and we show the
existence of a one-to-one correspondence between closed ideals of finite
codimension  and closed
submodules of finite codimension which preserves primary decompositions.
 At the end we consider the case of the Hardy space on the polydisc, and of
Bergman spaces on the unit ball and Reinhardt domains.
   
\vspace {15 mm}
 
In the development of multivariable operator theory there appears to
exist a fundamental connection between the study of invariant
subspaces for certain Hilbert spaces of holomorphic functions and the 
study of relatively closed ideals in rings of functions that are dense
in these spaces. Some results in this direction have been obtained
in [5] , [6] and [7]. Most of this has been done in the context of Hilbert
modules. In this paper we consider rings with topologies that 
carry some of the properties of those mentioned above and exhibit topological
properties for ideals in these rings. This general approach will give a
more transparent view of some results proved in [5].

\bigskip

{\centerline {\bf 1. Hilbert Nullstellensatz for closed ideals}}

\bigskip

Throughout the paper ${\cal R}$ will denote a commutative Noetherian
ring with unit, endowed with a topology $\tau $ for which  addition is
continuous and multiplication is separately continuous in each variable.

An example of such a ring is the ring ${\bf C}[z_1, z_2, \cdots ,z_n]$
of polynomials in $n$ variables with the topology induced by the Hardy
space of the polydisc $H^2({\bf D}^n)$, or the Bergman space $L^2_a(\Omega )$
of an open set $\Omega \subset {\bf C}^n$. Another example is the ring
${\cal O}(\stackrel {-}{{\bf B}})$ of analytic functions in a neighborhood
of the unit ball
  ${\bf B} \subset {\bf C}^n$ with the topology induced by $L^2_a(
{\bf B})$.

Let us remark that the ring is usually not complete in this topology.
In what follows we study properties of ideals that are closed
in the topology $\tau $.

\bigskip

{\bf Lemma 1.1. } The radical of a closed ideal is closed.

\medskip

{\bf Proof:} Let $I$ be closed, and $r(I)$ be its radical. By Proposition
7.14 in [2] there exists an integer $k$ such that $r(I)^k\subset I$.
Let $f_n$ be a sequence of elements in $r(I)$ converging to some $f$.
We want to prove that $f\in r(I)$. Since the multiplication is continuous 
in each variable, for every $g\in r(I)^{k-1}$, $f_ng\rightarrow fg$,
hence $fg\in I$, since $I$ is closed. This shows that in particular
$fg\in I$ for every $g\in r(I)^{k-1}$. Repeating the argument we get 
$ff_ng\rightarrow f^2g$ for any $g\in r(I)^{k-2}$, hence $f^2g\in I$ for
$g\in r(I)^{k-2}$.
Inductively we get $f^rg\in I$, for $g\in r(I)^{k-r}$ and $0\leq r\leq k$,   
so $f^k\in I$ which shows that $f\in r(I)$.

\bigskip
 
Since ${\cal R}$ is Noetherian, every ideal $I\subset {\cal R}$
has a (minimal) primary decomposition $I=Q_1\cap Q_2\cap \cdots \cap Q_m$
where each $Q_i$ is $P_i$-primary for some prime ideal $P_i$. The ideals $P_i$,
$1\leq i\leq m$ are called the prime ideals belonging to $I$ (or associated to
$I$).

\bigskip

{\bf Theorem 1.1.} If an ideal is closed, then every prime ideal belonging 
to it is closed.

\medskip

{\bf Proof:} If $I$ is closed and $f\in {\cal R}$, then the ideal
$(I:f)=\{g\in{\cal R}, gf\in I\}$ is closed. Indeed, if $g_n\in (I:f)$ and
$g_n\rightarrow g$ then $g_nf\rightarrow gf$, so $gf\in I$ which shows
that $g\in (I:f)$.

From Lemma 1.1. we get that $r((I:f))$ is closed for every $f$. By Theorem 4.5.
in [2] every prime ideal belonging to $I$ is of this form hence it is closed.

\bigskip

This gives a positive answer to a question raised in [5].

\medskip

{\bf Remark.} Since the closure of an ideal is an ideal, a maximal ideal
is either closed or dense.

\medskip

Given an ideal $J\subset {\cal R}$, the $J$-adic topology on ${\cal R}$ 
is the topology determined by the powers of $J$, so in this topology
the closure of a set $A\subset {\cal R}$ is $\cap _n(A+J^n)$. For more
details the reader can consult [2].

Let

\medskip

 ${\cal C}:=\{{\cal M}\subset {\cal R}, {\cal M} $ maximal ideal and
the ${\cal M}$-adic topology is weaker than $\tau \}$.

\medskip

We see that ${\cal C}$ consists of those maximal ideals ${\cal M}$ for which
${\cal M}^n$ is dense for every integer $n$.
The following result is a slightly modified version of Theorem 2.7 in [5].

\bigskip

{\bf Theorem 1.2.} If an ideal $I$ has the property that for every prime $P$ 
belonging to it there is ${\cal M}\in {\cal C}$ with $P\subset {\cal M}$,
then $I$ is closed.

\medskip

{\bf Proof:} Let $P_1,P_2, \cdots ,P_m$ be the prime ideals belonging to $I$,
$P_i\subset {\cal M}_i$, ${\cal M}_i\in {\cal C}$. If $J={\cal M}_1
{\cal M}_2\cdots {\cal M}_m$ then the $J$-adic topology is weaker than
$\tau$, so it suffices to prove that $I$ is closed in the $J$-adic 
topology. By Krull's Theorem ([2, Theorem 10.7]) and Proposition 4.7 in [2]
the latter is true if and only if $J+P_i\neq {\cal R}$, $1\leq i\leq m$.
This is clearly satisfied since $J+P_i\subset {\cal M}_i$, and the theorem is
proved.

\bigskip

{\bf Definition.} The pair $({\cal R}, \tau )$ is said to satisfy Hilbert's
Nullstellensatz if every ideal $I\subset {\cal R}$ is either dense,
or there exists ${\cal M}\in {\cal C}$ with $I\subset {\cal M}$.

\medskip

Let us remark that if $({\cal R},\tau )$ satisfies Hilbert's Nullstellensatz
then  any closed ideal is contained in a maximal closed ideal, which 
motivates the terminology.
By Krull's Theorem, $({\cal R}, \tau )$ satisfies Hilbert's Nullstellensatz
for every $J$-adic topology $\tau$. The ring ${\bf C}[z]$ with the topology
induced by $H^2({\bf D})$ or $L^2_a({\bf D})$ also satisfies this property.
 In [7], a class of strongly 
pseudoconvex domains $\Omega $
for which ${\cal O}(\stackrel{-}{\Omega})$ with the topology induced by
$L^2_a(\Omega )$ satisfies Hilbert's Nullstellensatz has been exhibited.

\bigskip

{\bf Lemma 1.2.} If $I$ and $J$ are two dense ideals in ${\cal R}$ then
$I\cdot J$ is dense.

\medskip

{\bf Proof:} Let $f_n\rightarrow 1$, $n\rightarrow \infty$, $f_n\in I$.
If $g\in J$ then $f_n g\rightarrow g$ which shows that $I\cdot J $ is dense
in $J$, hence in ${\cal R}$.

\bigskip

{\bf Theorem 1.3.} If $( {\cal R}, \tau )$ satisfies 
Hilbert's Nullstellensatz then an
ideal $I\subset {\cal R}$ is closed if and only if every prime belonging to
$I$ is closed. Moreover, the closure of an ideal in ${\cal R}$ is equal to the 
intersection of its primary components that are contained in closed maximal
ideals.

\medskip

{\bf Proof:} If $I$ is closed then every prime belonging to $I$ is closed by
Theorem 1.1. Conversely, if a prime belonging to $I$ is closed, then 
it is not dense, so it is included in an ${\cal M}\in {\cal C}$. The fact that
$I$ is closed now follows from Theorem 1.2.

For the second part, let $\stackrel{-}{I}$ be the closure of $I$ in 
${\cal R}$, and let $I=\cap_{i=1}^m Q_i$ be a (minimal) primary decomposition
of $I$, such that $Q_1,Q_2, \cdots ,Q_r$ are included in maximal ideals that
are in ${\cal C}$, hence closed, and $Q_{r+1},Q_{r+2},\cdots ,Q_m$ are not,
so they are dense. Then from the first part of the proof we get 
$\stackrel{-}{I}\subset \cap_{i=1}^nQ_i$.

On the other hand, by Lemma 1.2 the ideal $Q_{r+1}Q_{r+2}\cdots Q_m$ is dense
in ${\cal R}$, hence
$Q_{r+1}Q_{r+2}\cdots Q_m (Q_1\cap Q_2\cap \cdots \cap Q_r)$ is dense
in $Q_1\cap Q_2\cap \cdots \cap Q_r$, hence $I$ is dense in $ Q_1\cap Q_2 \cap
\cdots \cap Q_r$, so $\stackrel{-}{I}=Q_1\cap Q_2\cap \cdots \cap Q_r$.

\bigskip

In [5] and [6] the authors  state the following conjecture:

``Let ${\cal R}={\bf C}[z_1,z_2,\cdots ,z_n]$ be endowed with the topology
induced by $H^2({\bf D}^n)$. Then an ideal $I$ is closed if and only if 
every irreducible component of the zero set of $I$ intersects ${\bf D}^n$.''
  
The results above show that the conjecture is equivalent to the topological
Hilbert Nullstellensatz for the ring of polynomials with the topology
induced by the Hardy space.

\newpage

{\centerline {\bf 2. The case of the bidisk}}

\bigskip

In this section we prove the above mentioned
 conjecture  for the case of two variables.
Let us denote by ${\bf T}^2$ the 2-dimensional torus $\{(z_1,z_2)
\in {\bf C}^2; |z_i| =1\}$.

\bigskip

{\bf Lemma 2.1.} If $\alpha \in {\bf C}$ , $|\alpha |\geq 1$ and $1/2<r<1$
then for any $z$ with $|z|=1$ we have $|(z-\alpha )/(rz-\alpha )|\leq 2$.

\medskip

{\bf Proof:} The result follows from $|z-\alpha |\leq |z-\alpha /r|$.
The last inequality is obvious since in the triangle formed by
the points $z,\alpha $ and $\alpha /r$ the angle at $\alpha $ is
obtuse.

\bigskip

{\bf Lemma 2.2.} Let $p(z)=a_n(z-z_1)(z-z_2)\cdots (z-z_n)$ be such 
that $|z_i|\geq 1$. If $1/2<r<1$ then for any $z$ with $|z|
= 1$, $|p(z)/p(rz)|\leq 2^n$.

\medskip 

{\bf Proof:} The result follows by applying the previous lemma to
each of the factors in the decomposition of $p$.

\bigskip

{\bf Proposition 2.1.} Let $p\in {\bf C}[z_1,z_2]$ be a polynomial 
having no zeros in ${\bf D}^2$. Then $pH^2({\bf D}^2)$ is dense in
$H^2({\bf D}^2)$.	

\medskip

{\bf Proof:} Since every polynomial is a product of irreducible
polynomials, it suffices to prove the result for an irreducible polynomial.
Let us  show that 1 is in the closure of $pH^2({\bf D}^2)$.
 Without loss of generality we can assume that $p\not \in 
{\bf C}[z_2]$ (the case when $p$ is constant being obvious).
Write $p(z_1,z_2)=a_mz_1^m+\cdots +a_1z_1+a_0$ where $a_k\in {\bf C}[z_2]$.
If we fix $z_2^0$ on the circle, the polynomial in $z_1$, $p(z_1,z_2^0)$,
does not vanish identically; for otherwise it would be divisible by $z_2-
z_2^0$, contradicting the irreducibility. So, for each $z_2^0$ on the unit
 circle there is a maximal $k\leq m$ such that
 $a_k(z_2^0)\neq 0$.It follows that the polynomial $p$ 
has no zero $(z^0_1,z^0_2)$ with $|z^0_1|<1$ and $|z^0_2|=1$
since otherwise, by slightly perturbing $z_2$ to a point
 in ${\bf D}$ and using the 
continuous dependence of the roots on the coefficients, we could produce 
a zero for $p$ in ${\bf D}^2$. This shows   that if $0<r<1$ the polynomial
$p(rz_1,z_2)$ has no zeros in $\stackrel{-}{\bf D}^2$, thus $1/p(
rz_1,z_2)\in H^2({\bf D}^2)$.

Consider the sequence $f_k(z_1,z_2):=p(z_1,z_2)/
p((1-1/k)z_1,z_2)\in pH^2({\bf D}^2)$.
We want to show that $f_k\rightarrow 1$ in $H^2({\bf D}^2)$.
Fix $z_2^0$ on the unit
 circle and consider the  maximal $k\leq m$ such that
 $a_k(z_2^0)\neq 0$. By Lemma 2.2 we have
$|f_k(z_1,z_2^0)|\leq 2^k$ for any $z_1$ on the unit circle,
 hence $|f_k(z_1,z_2)|\leq 2^m$ for $(z_1,z_2)\in {\bf T}^2$.

Let us denote by $V(p)$ the zero set of $p$.
Since the 2-torus is not a complex algebraic variety, the set $V(p)\cap {\bf T}
^2$ has real dimension strictly less than 2, thus it has measure zero. 
For $\epsilon >0$, let us choose $W$ a neighborhood of this 
set, whose Lebesgue measure on the torus 
is equal to $\epsilon /(2(2^m+1)^2)$. Since $f_k\rightarrow 1$ uniformly
on ${\bf T}^n\backslash W$ we can choose $N\in {\bf N}$ such that 
for $k>N$, $\norm{f_k-1}_{2,{\bf T}^n\backslash W}^2<\epsilon  /2$.
 It follows  that for $k>N$, $\norm{f_k-1}_2^2<(2^m+1)^2 
\epsilon/(2(2^m+1)^2) +\epsilon /2=\epsilon$, which proves the assertion.

\bigskip

{\bf Remark.} We see that the only nontrivial situations when this 
result applies are those when the zero set of the polynomial
touches the boundary of ${\bf D}^n$. Here are some examples
of such polynomials: $z_1z_2-1$, $z_1+z_2-2$,
$2z_1z_2+z_1+z_2+2$.

\bigskip

{\bf Theorem 2.1.} The ring ${\bf C}[z_1,z_2]$, with the topology induced
by the Hardy space, satisfies the topological Hilbert Nullstellensatz.

\medskip

{\bf Proof:} Let us first prove the property for prime ideals.  
 Using the classical Hilbert Nullstellensatz we see that
the only maximal ideals in ${\cal C}$  are those corresponding to 
points in ${\bf D}^2$. Moreover, the other maximal ideals are dense. So
by Theorem 1.2 we only
 have to show that if $P$ is prime and $V(P)\cap {\bf D}^2=\emptyset $,
where $V(P)$ is the zero set of $P$,
then $P$ is dense in ${\bf C}[z_1,z_2]$. Standard results in dimension
theory (see [2]) show that $P$ is either maximal or principal. Indeed
the maximal length of a chain of nonzero prime ideals containing $P$ is 2.
If we have $P_0\subset P$ then $P$ is maximal and the result follows easily.
 If $P\subset P_1$ then $P$ must be principal since 
if $P$ is generated by $g_1,g_2, \cdots ,g_k$ we can take $g_1$
to be irreducible ( using the fact that $P$ is prime), and then 
the ideal generated by $g_1$ is included in $P$, so it must coincide 
with $P$. The density in the second case follows from Theorem 2.1.

If $I$ is an arbitrary ideal having no zeros in ${\bf D}^2$, let us show that
it is dense. If $P_1, P_2, \cdots ,P_n$ are the primes associated to it
then from what has been established above it follows that they are dense.
By Proposition 7.14 in [2] there exists an integer $k$
such that $(P_1P_2\cdots P_n)^k\subset I$. It follows from 
Lemma 1.2 that  $I$ itself
is dense in ${\bf C}[z_1,z_2]$, which proves the theorem.

\bigskip

As a direct consequence of this result and Theorem 1.3 we get

\medskip 

{\bf Corollary 2.1.} Let ${\bf C}[z_1,z_2]$ be endowed with the topology 
induced by the Hardy space. Then an ideal is closed if and only if
each of the irreducible components of its zero set intersects ${\bf D}^2$.

{\centerline {\bf 3. Ideals of finite codimension}}

\bigskip

Now let us suppose that ${\cal R}$ is also a $k$-algebra, where $k$ is an
algebraically closed field, and that the scalar multiplication is continuous
(not only separately continuous).  By the classical Hilbert
Nullstellensatz ([2, Corollary 5.24]), ${\cal R}/{\cal M}{\widetilde =}k$
for every maximal ideal ${\cal M}$. Let us also assume that the family
${\cal C}$ defined above consists of all closed maximal ideals. Thus in
this case a maximal ideal ${\cal M}$ is either dense, or the ${\cal M}$-adic
topology is weaker then the topology of ${\cal R}$.

\medskip

{\bf Example.} If we consider ${\bf C}[z_1,z_2,\cdots ,z_n]$ with
the topology induced by $H^2({\bf D}^n)$ then the condition above is
satisfied. In this case we have ${\cal C}={\bf D}^n$.

\medskip

 In a similar way as we proved Proposition 1.2 we can establish
the following result.

\bigskip

{\bf Lemma 3.1.} Given an ideal $I$ whose associated prime ideals are maximal,
$I$  is closed if and only if every
maximal ideal belonging to $I$ is closed.

\bigskip

{\bf Lemma 3.2.} An ideal $I\subset {\cal R}$ has  finite codimension
 in ${\cal R}$ if and only if every prime ideal belonging to $I$ is maximal.

\medskip

{\bf Proof:} Let $I=Q_1\cap Q_2\cap \cdots \cap Q_m$, $Q_i$
 ${\cal M}_i$-primary. By Proposition 7.14 in [2] there exists an integer $n$
such that 
${\cal M}_i^n\subset Q_i$, so $({\cal M}_1{\cal M}_2\cdots {\cal M}_m)
^n\subset I$. Since $dim{\cal R}/({\cal M}_1{\cal M}_2\cdots {\cal M}_m)^n<
\infty $, $I$ has finite codimension.

For the converse, let $P$ be a prime ideal belonging to $I$.
 Then $P$ has finite
codimension as well, so ${\cal R}/P$ is an integral domain that is finite over
$k$, and since $k$ is algebraically closed
 we must have ${\cal R}/P\stackrel{~}{=}k$, therefore $P$ is maximal.

\bigskip

Let $\widetilde {\cal R}$ be the closure of ${\cal R}$ in the topology $\tau$.
Since the multiplication is only separately continuous, $\widetilde {\cal R}$
is no longer a ring, but it is an ${\cal R}$-module. Each element
 $x\in {\cal R}$ induces a continuous multiplication morphism $T_x$ on
$\widetilde {\cal R}$. 
We shall denote by $\widetilde I$ the closure in $\widetilde
{\cal R}$ of an ideal $I$ in ${\cal R}$, to avoid confusion with 
$\stackrel{-}{I}$, the closure 
of $I$ in ${\cal R}$. Clearly $\widetilde I$ is a
closed submodule of $\widetilde {\cal R}$. Also, if $Y\subset 
\widetilde {\cal R}$ is a closed submodule, then $Y\cap {\cal R}$ is an
ideal that is closed in ${\cal R}$.
 
\medskip

{\bf Definition.}(see [2, page 58]) A submodule 
$Y \subset \widetilde {\cal R}$ is called primary in $\widetilde {\cal R}$
if $Y\neq \widetilde {\cal R}$ and every zero-divisor in $\widetilde {\cal R}
/Y$ is nilpotent.( An element $x\in {\cal R}$ is called zero-divisor if the
morphism induced by $T_x$ on $\widetilde {\cal R}/Y $ has nontrivial 
kernel, and nilpotent
if this morphism is nilpotent).

\medskip

 {\bf Remark. } If $Y\subset
 \widetilde {\cal R}$ is primary then $(Y:\widetilde
{\cal R})=\{x\in {\cal R}\ | \ T_x\widetilde {\cal R}\subset Y\}$ is primary,
so $P=r((Y:\widetilde {\cal R}))$ is prime. We say that $Y$ is $P$-primary.
Moreover, $(Y:\widetilde {\cal R})=Y\cap {\cal R}$ so $Y\cap {\cal R}$ is
also $P$-primary.

\medskip

Although every ideal in ${\cal R}$ has a primary decomposition, this is not 
true in general for the submodules of $\widetilde {\cal R}$. The next result
shows that this is true for closed submodules of $\widetilde {\cal R}$ 
of finite codimension.

\bigskip

{\bf Theorem 3.1.} There is a one-to-one correspondence between ideals in
  ${\cal R}$ whose associated prime ideals are maximal and closed in 
${\cal R}$, and  closed submodules in $\widetilde {\cal R}$ of finite
codimension, given by the maps $I\rightarrow \widetilde I$ and $Y\rightarrow
Y\cap {\cal R}$. Moreover, if $I=Q_1\cap Q_2\cap \cdots \cap Q_m$ is a 
( minimal) primary decomposition
 for $I$, then $\widetilde I=\widetilde Q_1 \cap
\widetilde Q_2 \cap \cdots \cap \widetilde Q_m$ is a ( minimal) primary
decomposition for $\widetilde I$.

\medskip

{\bf Proof: } By Lemmas 3.1 and 3.2 we have to show that the maps indicated 
above establish a one-to-one correspondence between ideals of finite 
codimension that are closed in the topology of ${\cal R}$ and closed submodules
of finite codimension in $\widetilde {\cal R}$.

If $Y$ is a closed submodule of $\widetilde {\cal R}$ of finite codimension
then ${\cal R}/(Y\cap {\cal R}){\widetilde =}\widetilde {\cal R}/Y$ since the
canonical map ${\cal R}\rightarrow \widetilde {\cal R}/Y$ is surjective,
${\cal R}$ being dense in $\widetilde {\cal R}$ and $\widetilde {\cal R}/Y$
being finite dimensional, and the kernel of this map is $Y\cap {\cal R}$.
On the other hand 
${\widetilde {\cal R}}/{\widetilde {(Y\cap {\cal R})}} 
{\widetilde =}{\cal R}/(Y\cap {\cal R})$,
hence $Y={\widetilde {(Y\cap {\cal R})}}$.
 Also for every ideal $I\subset {\cal R}$ that is 
closed in the topology of ${\cal R}$, $\widetilde I\cap {\cal R}=I$,
so the two maps are inverses of one another, and the one-to-one 
correspondence is proved.

Let $I=Q_1\cap Q_2 \cap \cdots \cap Q_m$ be a primary decomposition of $I$.
Then $\widetilde I\subset \widetilde Q_1\cap \widetilde Q_2 \cap \cdots \cap
\widetilde Q_m$, and since $\widetilde I\cap {\cal R}=\widetilde Q_1\cap
\widetilde Q_2\cap \cdots \cap \widetilde Q_m \cap {\cal R}=Q_1\cap Q_2\cap
\cdots Q_m=I$, by the first part of the proof the two must be equal.

In the commutative diagram below the horizontal arrows are isomorphisms

\medskip

{\centerline {${\cal R}/Q_i {\widetilde \rightarrow } {\widetilde {\cal R}}/
{\widetilde Q_i}$} }

\setlength{\unitlength}{1mm}
\begin{picture}(100,8)
\put(68,7){\vector(0,-1){8}}
\put(62,3){$T_x$}
\put(82,3){$T_x$}
\put(80,7){\vector(0,-1){8}}
\end{picture}

{\centerline {${\cal R}/Q_i {\widetilde \rightarrow } \widetilde {\cal R}/
\widetilde Q_i$}}

\medskip

\noindent so the fact that $Q_i$ is 
a primary ideal (hence a primary ${\cal R}$-module 
as well) implies that $\widetilde Q_i$ is a primary submodule of $\widetilde
{\cal R}$.

If the primary decomposition of $I$ is minimal let us show that 
the corresponding decomposition for
$\widetilde I$ is also minimal. Suppose that there exists $j$ such that
$\widetilde I=\widetilde Q_1\cap\cdots \widetilde Q_{j-1} \cap \cap 
\widetilde Q_{j+1}\cap \cdots \cap \widetilde Q_m$. Then $I=\widetilde I \cap
{\cal R}=Q_1\cap \cdots \cap Q_{j-1}\cap Q_{j+1}\cap \cdots \cap Q_m$,
which contradicts the minimality of the primary decomposition of $I$, and the
theorem is proved.

\bigskip

 From the proof it follows that the associated primes of $I$
and $\widetilde I$ coincide. Corollary 4.11 in [2] shows that in this case
the minimal primary decomposition is unique.

\medskip

{\bf Remark.}  In the case when the topology on ${\cal R}$ comes from a
norm, the first part of the theorem is contained in [5, Corollary 2.8].

\bigskip

{\centerline {\bf 4. Examples}}

In this last section we shall describe some examples to which Theorem 3.1
can be applied.

\medskip

1. For the case ${\cal R}={\bf C}[z_1, z_2, \cdots , z_n]$ and 
$\widetilde {\cal R} =H^2({\bf D}^n)$ the theorem has been proved by
Ahern and Clark in [1]. The primary closed
 ${\bf C}[z_1,z_2,\cdots ,z_n]$-submodules of finite codimension of $H^2(
{\bf D}^n)$ are those closed submodules $Y$ for which there exists a point 
$\lambda \in {\bf D}^n$ and a number $m\in {\bf N}$ such that $Y$ contains
the space of functions $f$ that satisfy 
$(\partial ^m/\partial z_1^m \partial ^m /\partial z_2^m \cdots
\partial ^m/\partial z_n^mf)(\lambda )=0$.  

\medskip

2. If ${\bf B} $ is the unit ball in ${\bf C}^n$, ${\cal R}=
{\cal O}(\stackrel{-}{{\bf B} })$ and $\widetilde {\cal R} =L^2_a({\bf B} )$
then ${\cal C}={\bf B} $, and if ${\cal M}$ is a maximal ideal corresponding
to a point in $\partial {\bf B} $ then ${\cal M}$ is dense in $L^2_a({\bf B} )$
(see [7]). This shows that the hypothesis of Theorem 2.1 is satisfied.
The primary closed ${\cal O}(\stackrel{-}{{\bf B} })$-modules of finite 
codimension have the same description as in the previous example.

\medskip

3. The last example concerns a certain class of Reinhardt domains.
 Let $0<p,q<\infty $. Following [4] we define 

\hspace{5mm} $\Omega _{p,q}=\{z\in {\bf C}^2 \ \ | \ \ |z_1|^p+|z_2|^q<1\}.$

We shall prove that ${\bf C}[z_1,z_2]$ with the topology induced by the
$L^2$-norm on $\Omega _{p,q}$ also satisfies the hypothesis of Theorem 3.1.

\bigskip

{\bf Lemma 4.1.} If $0\leq a \leq 1$ then the series 
\begin{eqnarray*}
\sum_{m,n\geq 0} \left(
\begin{array}{clcr}
m+n \\
m
\end{array}
\right)
a^m(1-a)^n
\end{eqnarray*}
diverges.

\medskip

{\bf Proof:} By symmetry, we can assume $a< 1$. We have 
\begin{eqnarray*}
\sum_{m,n\geq 0}\left(
\begin{array}{clcr}
m+n \\
m
\end{array}
\right)
a^m(1-a)^n=\sum_{p\geq 0}\sum_{q\leq p}\left(
\begin{array}{clcr}
p \\
q
\end{array}
\right)
a^q(1-a)^{p-q}=
\end{eqnarray*}
\begin{eqnarray*}
=\sum_{p\geq 0}(1-a)^p\sum_{q\leq p}\left(
\begin{array}{clcr}
p \\
q
\end{array}
\right)( \frac{a}{1-a})
^q=\sum_{p\geq 0}(1-a)^p\cdot \frac{1}{(1-a)^p}=\infty
\end{eqnarray*}

\bigskip

{\bf Corollary 4.1.} If $0\leq a \leq 1$ then the series 
\begin{eqnarray*}
\sum_{m,n> 0}\left(
\begin{array}{clcr}
2(m+n)-2 \\
2m-1
\end{array}
\right )
a^{2m-1}(1-a)^{2n-1}
\end{eqnarray*}
\noindent diverges.

\medskip

Let $B$ denote the beta function, defined by
 $B(r,s):=\Gamma (r)\Gamma (s)/\Gamma (r+s)$, $r,s>0$.

\noindent By [3],
\begin{eqnarray*}
B(r,s)=\int_0^1 t^{r-1}(1-t)^{s-1}dt.
\end{eqnarray*}

\bigskip

{\bf Lemma 4.2.} If $(z_1,z_2)\in \partial \Omega _{p,q}$ then the series
\begin{eqnarray*}
\sum_{r_1,r_2 \geq 0}
\frac{|z_1|^{2r_1}|z_2|^{2r_2}}{B({\frac{2r_1+2}{p}},{\frac{2r_2+2}{q}}+2)}
\end{eqnarray*}

\noindent diverges. 

\medskip

{\bf Proof:} Since $B(r,s)$ is a decreasing function in $r$ and $s$, by
taking the subseries corresponding to indices with the property that
$pm\leq r_1\leq pm+1$ and $qn\leq r_2 \leq qn+1$ we get the inequality
\begin{eqnarray*}
\sum _{r_1, r_2\geq 0}
\frac{|z_1|^{2r_1}|z_2|^{2r_2}}{B({\frac{2r_1+2}{p}},{\frac{2r_2+2}{q}}+2)}
\geq |z_1|^{p+2}|z_2|^{q+2}\sum_{m,n\geq 0}
{\frac{|z_1|^{(2m-1)p}|z_2|^{(2n-1)q}}{B(2m,2n)}}.
\end{eqnarray*}

\noindent Since $|z_1|^q=1-|z_2|^p$ the latter series is equal to
 $\sum_{m,n\geq 0}(|z_1|^p)^{2m-1}(1-|z_1|p)^{2n-1}/B(2m,2n)$.
 Applying Corollary
4.1 for $a=|z_1|^p$, and using the fact that 
$B(2m,2n)=(2m-1)!(2n-1)!/(2m+2n-1)!<1/(^{2(m+n)-2}_{\  2m-1} ) $ (see [3])
 we get the desired result.

\bigskip

{\bf Proposition 4.1.} If we endow the ring ${\bf C}[z_1,z_2]$
with the topology induced by  $L^2(\Omega _{p,q})$, then a
maximal ideal ${\cal M}\subset {\bf C}[z_1,z_2]$ is either
dense, or the ${\cal M}$-adic topology is weaker then the topology of the ring.

\medskip

{\bf Proof:} By [4, Theorem 4.9 c) and  Example 5.2], an ideal 
${\cal M}=(z_1-w_1,z_2-w_2)$ is dense in ${\bf C}[z_1,z_2]$ 
with the $L^2$-norm if and only if 
\begin{eqnarray*}
{\frac{2(\pi)^2}{p}} \sum_{r_1,r_2\geq 0}(r_2+1)
\frac{|w_1|^{2r_1}|w_2|^{2r_2}}{B({\frac{2r_1+2}{p}},{\frac{2r_2+2}{q}}+1)}.
\end{eqnarray*}
\noindent diverges.

By Lemma 4.2 this series diverges if $(w_1,w_2)\in
\partial \Omega _{p,q}$, so if $w_1, w_2 \not \in \Omega _{p,q}$
the ideal ${\cal M}$ is dense.

It is not difficult to check that ${\cal M}^m$ is closed whenever 
$(w_1,w_2)\in \Omega _{p,q}$ and $m\in {\bf N}$. Therefore 
${\cal C}=\Omega _{p,q}$, and if ${\cal M}\not \in {\cal C}$ then ${\cal M}$ 
is dense.
 
\vspace{8mm}
 
{\bf Bibliography:}

\medskip

[1] Ahern, P., Clark, D., N., {\em Invariant subspaces and analytic 
continuation in several variables}, J. Math. Mech., 19(1969/1970), 963--969.

[2] Atiyah, M., F., Macdonald, I., G., {\em Introduction to commutative
 algebra}, Addison-Wesley, 1969.

[3] Carlson, B., {\em Special functions 
of applied mathematics}, Academic Press,
New York, 1977.

[4] Curto, R., Salinas, N., {\em Spectral properties of cyclic subnormal
m-tuples}, Amer. J. Math., Vol. 107, 1(1985), 113--138.
 
[5] Douglas, R., Paulsen, V., Sah, C., Yan, K., {\em Algebraic reduction 
and rigidity for Hilbert modules}, Amer. J. Math., to appear.

[6] Paulsen, V., {\em Rigidity theorems in spaces of analytic functions}, Proc.
Symp. Pure Math., Vol. 51(1991), Part 1.

[7] Putinar, M., Salinas, N., {\em Analytic transversality and Nullstellensatz
in Bergman space}, Contemporary Math., Vol. 137, 1992.

\vspace{10mm}

\noindent Department of Mathematics, 
The University of Iowa, Iowa City, IA 52242

\noindent {\em E-mail: rgelca@math.uiowa.edu}

\medskip

and

\medskip

\noindent Institute of Mathematics of the
 Romanian Academy, P.O.Box 1-764, Bucharest
70700, Romania.

\end{document}